\documentclass[reqno]{amsart}
\usepackage{mathtools}
\usepackage{dsfont}
\usepackage{amsmath,amssymb,amsthm,amsthm}
\usepackage{url,setspace,amsfonts,amstext,amsopn}
\usepackage{graphicx}
\usepackage{color}
\usepackage{enumerate}
\usepackage[utf8]{inputenc}

\usepackage{float}

\newtheorem{theorem}{Theorem}
\newtheorem*{theorem*}{Theorem}
\newtheorem{lemma}{Lemma}
\newtheorem{proposition}{Proposition}
\newtheorem{corollary}{Corollary}
\theoremstyle{remark}
\newtheorem{remark}{Remark}

\usepackage[colorlinks]{hyperref} 
\hypersetup{
    colorlinks=salmon,
    citecolor=blue, 
   linkcolor=yellowgreen, 
 linkcolor=TRKY,
}

\definecolor{mahogany}{cmyk}{0, 0.77, 0.87, 0}
\definecolor{salmon}{cmyk}{0, 0.53, 0.38, 0}
\definecolor{melon}{cmyk}{0, 0.46, 0.50, 0}
\definecolor{yellowgreen}{cmyk}{0.44, 0, 0.74, 0}
\definecolor{brickred}{cmyk}{0, 0.89, 0.94, 0.28}
\definecolor{OliveGreen}{cmyk}{0.64, 0, 0.95, 0.40}
\definecolor{RawSienna}{cmyk}{0, 0.72, 1.0, 0.45}
\definecolor{ZurichRed}{rgb}{1, 0, 0} 
\definecolor{TRKY}{rgb}{1, 0, 0.5}
\definecolor{blue}{rgb}{0, 0, 1}

\title
[Refined Eigenvalue Bounds on the Dirichlet Fractional Laplacian]
{Refined Eigenvalue Bounds on the Dirichlet Fractional Laplacian}

\author{Selma  Y{\i}ld{\i}r{\i}m Yolcu}
\address
[Selma  Y{\i}ld{\i}r{\i}m  Yolcu]
{Bradley University, Department of Mathematics, Peoria, IL 61625 USA}
\email{syildirimyolcu@bradley.edu}

\author{T\"{u}rkay Yolcu}
\address
[T\"{u}rkay Yolcu]
{Bradley University, Department of Mathematics, Peoria, IL 61625 USA}
\email{tyolcu@bradley.edu}

\date{\today}

\begin{document}
\nocite{*}

\begin{abstract} The purpose of this article is to establish new lower bounds for the sums of powers of eigenvalues of the 
Dirichlet fractional Laplacian operator 
$(-\Delta)^{\alpha/2}|_{\Omega}$ restricted to a bounded domain $\Omega\subset{\mathbb R}^d$ with $d=2,$ $1\leq \alpha\leq 2$ and $d\geq 3,$ $0< \alpha\le 2$. 
Our main result yields a sharper lower bound, in the sense of Weyl asymptotics, for the Berezin-Li-Yau type inequality improving the previous result in \cite{ST6}. Furthermore, we give a result improving the bounds for analogous elliptic operators in \cite{Kim}. 
\end{abstract}

\maketitle

\section{introduction}

Fractional Laplacian operators are usually considered as the prototype of the non-local operators \cite{ChenKimSong} and have recently garnered much attention in many applications in mathematics and physics. Related problems usually lie at the interface of probability, stochastic processes, partial differential equations and spectral theory such as \cite{BanYol, ChenKimSong, ChenSong}. For some applications, we refer the reader to graphene models \cite{HajjMeh}, obstacle problems \cite{Silv}, non-local minimal surfaces \cite{CafRoSav}.

In this article, we study eigenvalues of the fractional Laplacian operator $(-\Delta)^{\alpha/2}$ defined by
\begin{equation}\label{cpp}
\begin{split}
(-\Delta)^{\alpha/2}\,\phi_j &=\lambda_j^{(\alpha)}\,\phi_j \quad \hbox{in} \,\,\Omega,\\ 
 \phi_j  &=0 \quad\hbox{in} \,\,{\mathbb R}^d\backslash \Omega 
\end{split}
\end{equation}
where $\Omega$ is a bounded connected domain with smooth boundary in ${\mathbb R}^d$, for either $d=2,$ $1\le \alpha\le 2$ or $d\ge 3$ and $0< \alpha \le 2.$ Since $\Omega$ is bounded, the spectrum of the fractional Laplacian is discrete and eigenvalues $\{\lambda_j^{(\alpha)}\}_{j=1}^{\infty}$ (including multiplicities) can be sorted in an increasing order.

Unlike Laplacian, fractional Laplacian is a non-local operator such that for suitable test functions, including all functions in $u\in C_0^{\infty}({\mathbb R}^d)$, it is defined as 
\begin{equation}\label{fracdef1}
(-\Delta)^{\alpha/2}u(x)={\mathcal A}_{d,\alpha}\lim_{\epsilon\to0^+}\int_{\{|y|>\epsilon\}}\frac{u(x+y)-u(x)}{|y|^{d+\alpha}}dy,
\end{equation}
where ${\mathcal A}_{d,\alpha}$ is a well-known positive normalizing constant. 
From a probabilistic point of view, $(-\Delta)^{\alpha/2}$ can be considered as the infinitesimal generator of the semigroup of the symmetric $\alpha-$stable process, denoted by $X_t$, with the characteristic function
\begin{equation}\label{probdefn}
e^{-t|\mu|^{\alpha}}=E(e^{i\mu\cdot X_t})= \int_{{\mathbb R}^d}e^{i\mu\cdot y}p_t^{(\alpha)}(y)dy,\qquad t>0, \quad \mu\in{\mathbb R}^d,
\end{equation}
where  $p_t^{(\alpha)}(x,y)=p_t^{(\alpha)}(x-y)$ is called the transition density of the stable process (or the heat kernel of the fractional Laplacian). While explicit formulae for the transition density of  symmetric $\alpha$-stable processes are only available for  the Cauchy process ($\alpha=1$) and the Brownian motion ($\alpha=2$), these processes share many of the basic properties of the Brownian motion.  Another process of importance is the Holtsmark distribution ($\alpha=3/2$) which is used to model gravitational fields of stars (See e.g., \cite{Zolo}). Stable processes do not have continuous paths which is related to non-locality of the fractional Laplacian operator \cite{BanYol,Bogetal}.

When fractional Laplacians involved, some of the known methods fail because of the fractional power and non-locality of such operators. This drawback can be evaded by using the Fourier transform definition. Recall that the Fourier transform and its inverse are defined as
   $${\mathcal F}[u](\mu) =\hat{u}(\mu)=c_d\int_{{\mathbb R}^d}e^{-i\mu\cdot z}\,u(z)\,dz, \quad{\mathcal F}^{-1}[u](z)=c_d\int_{{\mathbb R}^d}e^{i\mu\cdot z}\,u(\mu)\,d\mu,$$
 where $c_d=(2\pi)^{-\frac{d}{2}}$ is the normalizing constant.  Interestingly, the fractional Laplacian operator on $\Omega\subset{\mathbb R}^d$ can be defined as a pseudo-differential operator as
 \begin{equation}\label{defnbilapfourier}
 (-\Delta)^{\alpha/2}|_{\Omega}u = {\mathcal F}^{-1}\left[|\mu|^{\alpha }{\mathcal F}[u]\right],\qquad 0<\alpha \le2, \quad u\in H_0^{\alpha/2}(\Omega).
 \end{equation}
 Here, $H_0^{\alpha/2}(\Omega)$ denotes the Sobolev
 space of order $\alpha/2$. When $\Omega={\mathbb R}^d$, one can look at Proposition 3.3. \cite{Vald} for the proof of the equivalence between the definitions in \eqref{fracdef1} and \eqref{defnbilapfourier}.
 
There is an extensive literature devoted to the inequalities involving the eigenvalues of the Dirichlet Laplacian operator, which can be regarded as the fractional Laplacian when $\alpha=2$. One may consult the articles \cite{AshBen, AshBenLaug, AshLaug, FraGei, Henrot, KovVuWei, LapGeiWei, LapWei, Weidl} and references therein for a through literature review. It is worth pausing here for a moment to consider the Dirichlet Laplacian results relevant to our main result. The first such result is the Li-Yau inequality that provides a lower bound for the sums of eigenvalues, sharp in the sense of Weyl asymptotics \cite{LiYau}. The authors generalized this result in \cite{ST2} by obtaining the following Li-Yau type inequality for the eigenvalues of the Diriclet fractional Laplacian operator:
 \begin{equation}
\sum_{j=1}^k\lambda_j^{(\alpha)}
\ge (4\pi)^{\frac{\alpha}{2}}\,\frac{d}{ \alpha+d }
\left(\frac{\Gamma\left(1+\frac{d}{2}\right)}{ |\Omega|}\right)^{\frac{\alpha }{d}}k^{1+\frac{\alpha }{d}},\label{BLY}
\end{equation}
where $|\Omega|$ represents the volume of $\Omega$ and $\Gamma(x)$ denotes the Gamma function $\Gamma(x)=\int_0^{\infty} t^{x-1}e^{-t}\,dt$ for $x>0.$
One may also look at \cite{HarYil} for a Li-Yau type inequality involving the eigenvalues of the (massless) Klein--Gordon square root operators $(-\Delta)^{1/2}|_{\Omega}$, (i.e., the case $\alpha=1$). To look at this inequality from a different perspective, one may 
take the Legendre transform of the following result by Laptev \cite{Laptev2} and obtain \eqref{BLY}: 
\begin{equation}
\sum_{j}(z-\lambda_j^{(\alpha)})_+\leq (4\pi)^{-\frac{d}{2}} \frac{\alpha}{\alpha+d}\frac{|\Omega|}
{\Gamma\left(1+\frac{d}{2}\right)} z^{1+\frac{d}{\alpha}}.\label{AL}
\end{equation}
When we set $\alpha=2$ in \eqref{AL}, we recover an earlier result by Berezin \cite{Bere}, which supplies the Li-Yau inequality after an application of the Legendre transform. Thus, in what follows, we call \eqref{BLY} as the Berezin-Li-Yau inequality. 

In \cite{Melas}, Melas proved the following improvement to the Berezin-Li-Yau inequality ($\alpha=2$):
\begin{equation}\label{melasimp}
\sum_{j=1}^k\lambda_j^{(2)}
\ge 4\pi\,\frac{d}{ 2+d }
\left(\frac{\Gamma\left(1+\frac{d}{2}\right)}{ |\Omega|}\right)^{\frac{2 }{d}}k^{1+\frac{2}{d}}+\frac{1}
 {24(2+d )}\frac{|\Omega|}{{\mathcal I}(\Omega)}\,k,
\end{equation}
where ${\mathcal I}(\Omega)$, the moment
 of inertia, is defined by
 $${\mathcal I}(\Omega)=\min_{y\in{\mathbb R}^d}\int_{\Omega} |z-y|^2\,dz.$$
 By a translation of the origin and a rotation of axes if necessary, in the sequel, we assume that the origin is the center of mass of  $\Omega$ and that
 \begin{equation}
 {\mathcal I}(\Omega)=\int_{\Omega} |z|^2\,dz.\label{ID}
 \end{equation}
Melas type bounds and their many variants and extensions have recently attracted a lot of attention, see for instance \cite{Ilyin, KovWei, Selma, ST4, ST2,  ST5, ST3, ST6, Turkay}. In particular, in \cite{ST2}, the authors obtained a refinement of \eqref{BLY}, stating that
\begin{equation}
\begin{split}\sum_{j=1}^k\lambda_j^{(\alpha)}
& \ge (4\pi)^{\frac{\alpha}{2}}\,\frac{d}{ \alpha+d }
\left(\frac{\Gamma\left(1+\frac{d}{2}\right)}{ |\Omega|}\right)^{\frac{\alpha }{d}}k^{1+\frac{\alpha }{d}}
\\ &\quad +
\frac{\alpha}
 {48( \alpha+d )}\frac{|\Omega|^{1-\frac{\alpha-2}{d}}\,\Gamma\left(1+\frac{d}{2}\right)^{\frac{\alpha-2}{d}}}{(4\pi)^{1-\frac{\alpha}{2}}\,{\mathcal I}(\Omega)}\,k^{1+\frac{\alpha-2}{d}}.
 \end{split}\label{ccmbiz}
 \end{equation}
 Remark that $\alpha=2$ in \eqref{ccmbiz} recovers Melas' bound  in \cite{Melas} for the Dirichlet Laplacian eigenvalues. Unfortunately, it is not easy to take Legendre transform of \eqref{ccmbiz} to find a similar improved Berezin type bound in the case of the fractional Laplacian.

 For $0<\ell\leq 1$, authors also proved in \cite{ST6} that
 \begin{equation}\label{imbound2e}
 \begin{split}
 \sum_{j=1}^k \left(\lambda_j^{(\alpha)}\right)^\ell
 &\ge (4\pi)^{\frac{\alpha\ell}{2}}\,\frac{d}{ \alpha\ell+d }
 \left(\frac{\Gamma\left(1+\frac{d}{2}\right)}{ |\Omega|}\right)^{\frac{\alpha\ell}{d}} k^{1+\frac{\alpha\ell }{d}} \\
 &\quad+\frac{\alpha\ell}
  {16(\alpha\ell+d )} \frac{|\Omega|^{1-\frac{\alpha\ell-2}{d}}\,\Gamma\left(1+\frac{d}{2}\right)^{\frac{\alpha\ell-2}{d}}}{(4\pi)^{1-\frac{\alpha\ell}{2}}\,{\mathcal I}(\Omega)}\,k^{1+\frac{\alpha\ell-2}{d}}
 \\& \quad - \frac{\alpha\ell}{640(\alpha\ell+d )}
 \frac{|\Omega|^{2-\frac{\alpha\ell -4}{d}}\,\Gamma\left(1+\frac{d}{2}\right)^{\frac{\alpha\ell -4}{d}}}{ {(4\pi)^{2-\frac{\alpha\ell}{2}}\,\mathcal I}(\Omega)^{2}}\,k^{1+\frac{\alpha\ell -4}{d}}. 
 \end{split}
 \end{equation}
The main purpose of this article is to obtain analogous but sharper bounds with Dirichlet Laplacian replaced by Dirichlet fractional Laplacian. The first step, inspired by results in \cite{Selma, ST4, ST2, ST6, Turkay}, is establishing the following result:
\begin{proposition}\label{impBLY} 
For $k\ge 1$, and either $1\le\alpha\le 2$ and $d=2$ or  $0<\alpha\le 2$ and $d\ge 3,$ the eigenvalues $\{\lambda_j^{(\alpha)}\}_{j=1}^{\infty}$ of the fractional Laplacian operator \eqref{cpp} defined on $\Omega \subset {\mathbb R}^d$ satisfy 
 \begin{equation}\label{imbound1}
 \begin{split}
 \sum_{j=1}^k\lambda_j^{(\alpha)}
 &\ge(4\pi)^{\frac{\alpha}{2}}\,\frac{d}{ \alpha+d }
 \left(\frac{\Gamma\left(1+\frac{d}{2}\right)}{ |\Omega|}\right)^{\frac{\alpha }{d}}k^{1+\frac{\alpha }{d}}
 \\& \quad+\frac{\alpha}
  {2( \alpha+d)} \frac{|\Omega|^{\frac{1}{2}-\frac{\alpha-1}{d}}\,\Gamma\left(1+\frac{d}{2}\right)^{\frac{\alpha-1}{d}}}{(4\pi)^{\frac{1}{2}-\frac{\alpha}{2}}\,{\mathcal I}(\Omega)^{\frac{1}{2}}}\,k^{1+\frac{\alpha-1}{d}}
 \\& \quad - \frac{5\alpha}
  {16( \alpha+d )}\frac{|\Omega|^{1-\frac{\alpha-2}{d}}\,\Gamma\left(1+\frac{d}{2}\right)^{\frac{\alpha-2}{d}}}{(4\pi)^{1-\frac{\alpha}{2}}\,{\mathcal I}(\Omega)}\,k^{1+\frac{\alpha-2}{d}}
  \\& \quad+\frac{\alpha}{16( \alpha+d )}
 \frac{|\Omega|^{\frac{3}{2}-\frac{\alpha -3}{d}}\,\Gamma\left(1+\frac{d}{2}\right)^{\frac{\alpha -3}{d}}}{ {(4\pi)^{\frac{3}{2}-\frac{\alpha}{2}}\,\mathcal I}(\Omega)^{\frac{3}{2}}}\,k^{1+\frac{\alpha -3}{d}}. 
 \end{split}
 \end{equation}
\end{proposition}
Note that the constants in the leading terms on the right side of \eqref{imbound1}, which is a fractional version of Weyl's law, are optimal due to a classical result of Blumenthal and Getoor \cite{BluGet}. While the constant in the second term in \eqref{imbound1} is still sub-optimal, the estimate in \eqref{imbound1} is substantially stronger than previous known results in \cite{ST2,ST6,Turkay}. In addition, we recover the lower bounds in the case of the Dirichlet Laplacian \cite{ST4} when we set $\alpha=2$ in \eqref{imbound1}. 

In \cite{ChenSong}, Chen and Song obtained that
\begin{equation}
\lambda_j^{(\alpha\ell)}\le \left(\lambda_j^{(\alpha)}\right)^\ell\label{CS}
\end{equation}
for each $j$ and any constant $\ell\in(0,1]$. Thus, Proposition \ref{impBLY} along with an application of \eqref{CS} leads to our principal result:
\begin{theorem}\label{impBLYT} 
For $k\ge 1$, $0<\ell\le 1$ and either $1\le\alpha\le 2$ and $d=2$ or  $0<\alpha\le 2$ and $d\ge 3,$ the eigenvalues $\{\lambda_j^{(\alpha)}\}_{j=1}^{\infty}$ of the fractional Laplacian operator \eqref{cpp} defined on $\Omega \subset {\mathbb R}^d$ satisfy 
 \begin{equation}\label{imbound21}
 \begin{split}
 \sum_{j=1}^k\left(\lambda_j^{(\alpha)}\right)^{\ell}
 &\ge(4\pi)^{\frac{\alpha\ell}{2}}\,\frac{d}{ \alpha\ell+d }
 \left(\frac{\Gamma\left(1+\frac{d}{2}\right)}{ |\Omega|}\right)^{\frac{\alpha\ell}{d}}k^{1+\frac{\alpha\ell }{d}}
 \\& \quad+\frac{\alpha\ell}
  {2( \alpha\ell+d)} \frac{|\Omega|^{\frac{1}{2}-\frac{\alpha\ell-1}{d}}\,\Gamma\left(1+\frac{d}{2}\right)^{\frac{\alpha\ell-1}{d}}}{(4\pi)^{\frac{1}{2}-\frac{\alpha\ell}{2}}\,{\mathcal I}(\Omega)^{\frac{1}{2}}}\,k^{1+\frac{\alpha\ell-1}{d}}
 \\& \quad - \frac{5\alpha\ell}
  {16( \alpha\ell+d )}\frac{|\Omega|^{1-\frac{\alpha\ell-2}{d}}\,\Gamma\left(1+\frac{d}{2}\right)^{\frac{\alpha\ell-2}{d}}}{(4\pi)^{1-\frac{\alpha\ell}{2}}\,{\mathcal I}(\Omega)}\,k^{1+\frac{\alpha\ell-2}{d}}
  \\& \quad+\frac{\alpha\ell}{16( \alpha\ell+d )}
 \frac{|\Omega|^{\frac{3}{2}-\frac{\alpha\ell -3}{d}}\,\Gamma\left(1+\frac{d}{2}\right)^{\frac{\alpha\ell -3}{d}}}{ {(4\pi)^{\frac{3}{2}-\frac{\alpha\ell}{2}}\,\mathcal I}(\Omega)^{\frac{3}{2}}}\,k^{1+\frac{\alpha\ell -3}{d}}. 
 \end{split}
 \end{equation}
\end{theorem}

On a side note, the proof of Proposition \ref{impBLY} consists of a very delicate application of Steffensen's type inequalities, which is mainly about the comparison of the integrals on the subsets of interval $[0,\infty)$.

The article is structured as follows: In Section \ref{sec:BLY2}, we revisit the relevant facts about the eigenvalues and eigenfunctions of the fractional Laplacian operator. After providing the proof of an auxilliary lemma that plays a crucial role in proving Proposition \ref{impBLY}, we present the proof of our main results in Section \ref{sec:impBLY}. Finally, we end Section \ref{sec:impBLY} with a remark which extends the main result even further for some elliptic operators studied in \cite{Kim,SV}. Please see Remark \ref{rmkgen} for details.
\section{Preliminaries}\label{sec:BLY2}
In this section, we give an overview of the definitions and tools that are essential to establish the estimates in \eqref{imbound1}. Even though, these were previously studied in \cite{Selma, ST2, ST6, Turkay}, we include them so that the article is self-contained. 
By using Plancherel's theorem, one can show that the set of Fourier transforms $\{\hat{\phi}_j\}_{j=1}^{\infty}$ of $\{\phi_j\}_{j=1}^{\infty}$ forms an orthonormal set
in $L^2({\mathbb R}^d)$ since the set of eigenfunctions $\{\phi_j\}_{j=1}^{\infty}$ is an orthonormal set in $L^2(\Omega).$ To ease the notation, we set
\begin{equation}\label{Uk}
\Phi_k(\mu):=\sum_{j=1}^k|\hat{\phi}_j(\mu)|^2=\frac{1}{(2\pi)^{d}}\sum_{j=1}^k\left|\int_{\Omega} 
e^{-i z\cdot\mu}\phi_j(z)\,dz\right|^2 \ge 0.
\end{equation}
Because the support of $\phi_j$ is $\Omega,$ the integral is taken over $\Omega$ instead of $\mathbb{R}^d$. Interchanging the sum and integral and using $\|\hat{\phi}_j\|_2=1,$  we derive
\begin{equation}\label{Fintk}
\int_{{\mathbb R}^d}\Phi_k(\mu)\, d\mu=k.
\end{equation}
Observe that
\begin{equation}\label{Uksumlambda}
\begin{split}
\sum_{j=1}^k \lambda_j^{(\alpha)}&=\sum_{j=1}^k \langle \phi_j,\lambda_j^{(\alpha)} \phi_j\rangle   =\sum_{j=1}^k \langle \phi_j,(-\Delta)^{\alpha/2}\phi_j\rangle \\ & =\sum_{j=1}^k \langle \phi_j,\mathcal{F}^{-1}[|\mu|^{\alpha }\mathcal{F}[\phi_j]]\rangle =\sum_{j=1}^k \int_{{\mathbb R}^d}|\mu|^{\alpha }\,|\hat{\phi}_j(\mu)|^2\,d\mu \\ & =\int_{\mathbb{R}^d} |\mu|^{\alpha }\, \Phi_k(\mu)\, d\mu.
\end{split}
\end{equation}
Application of Bessel's inequality gives an upper bound for $\Phi_k$:
\begin{equation}
\Phi_k(\mu)\le \frac{1}{(2\pi)^{d}}\int_{\Omega}\left|e^{-iz\cdot\mu}\right|^2\,dz = \frac{|\Omega|}{(2\pi)^{d}}:=\Omega_d.\label{fksiless}
\end{equation}
Next, we find an estimate for $|\nabla \Phi_k|$.  Note that
\begin{equation}
\sum_{j=1}^{k}|\nabla \hat{\phi}_j(\mu)|^2\leq\frac{1}{(2\pi)^{d}}\int_{\Omega}\left|iz e^{-iz\cdot\mu}\right|^2\,dz
=\frac{{\mathcal I}(\Omega)}{(2\pi)^{d}}.\label{bugrad}
\end{equation}
In view of H{\" o}lder's inequality and utilizing \eqref{fksiless} and \eqref{bugrad}, we arrive at the following uniform bound:
\begin{equation}\label{gradUk}
\begin{split}
|\nabla \Phi_k(\mu)|&\le 2\left(\sum_{j=1}^{k}|\hat{\phi}_j(\mu)|^2\right)^{1/2}
\left(\sum_{j=1}^{k}|\nabla \hat{\phi}_j(\mu)|^2\right)^{1/2}\\ & \le 2(2\pi)^{-d}\sqrt{|\Omega|\, {\mathcal I}(\Omega)}:=\beta.
\end{split}
\end{equation}
Let $B_R(z):=\{y\in\mathbb{R}^d : |y-z|\le R\}$ designate the ball of radius $R$ centered at $z$ in $\mathbb{R}^d$ and $\omega_d$ denotes the volume of $d$ dimensional unit ball $B_1(z)$ in $\mathbb{R}^d$ given by
\begin{equation}
\omega_d=\dfrac{\pi^{\frac{d}{2}}}{\Gamma\left(1+\frac{d}{2}\right)}.\label{wd}
\end{equation}
 Now assume that $R$ is such that $|\Omega|=\omega_dR^d$. That is, $B_R(0)$ is the symmetric rearrangement of $\Omega.$ 
Note that
\begin{equation*}
{\mathcal I}(\Omega)\geq \int_{B_R(0)}|z|^2\,dz=\frac{d\omega_d}{d+2}R^{d+2}
=\frac{d}{d+2}\omega_d^{-\frac{2}{d}}|\Omega|^{\frac{d+2}{d}},\label{IDineq}
\end{equation*}
roughly giving
\begin{equation} \beta
\geq (2\pi)^{-d}\,\omega_d^{-\frac{1}{d}}\,|\Omega|^{\frac{d+1}{d}}.\label{boundform}
\end{equation}
Let $\Phi_k^*(\mu)$ denote the decreasing radial rearrangement of $\Phi_k(\mu).$
There exists a real valued absolutely continuous function $\varphi_k:[0,\infty)\to [0,\Omega_d]$ such that 
\begin{equation}\Phi_k^*(\mu)=\varphi_k(|\mu|).\label{zetak}
\end{equation} 
By using P{\'o}lya-Szeg{\"o} inequality, one can show that
\begin{equation}0\le -\varphi_k'(t)\le \beta.\label{bdv}
\end{equation}
For more details, see for example \cite{ST2}.

\section{Proof of Proposition \ref{impBLY}} \label{sec:impBLY}
Before diving into the proof of the main results, we present the following surprising sharper inequality which will be the main ingredient in the proof of the sharper lower bound in \eqref{imbound1}. Our method of proof has been previously explored in several articles, for instance \cite{Melas,Selma,ST4, ST5, ST6}, with crucial differences.
\begin{lemma}\label{lem:ineq} 
For either $d=2$ and $1\le\alpha\le 2$ or $3\le d\in \mathbb{N}$ and $0<\alpha\le 2,$  $a>0,$ $b>0,$ we have the following inequality
\begin{equation}
a^{d+\alpha}\ge \frac{d+\alpha}{d}a^{d} b^{\alpha}-\frac{\alpha}{d}b^{d+\alpha}+\frac{\alpha}{d} \,b^{d+\alpha-3}\,\left(2a+b\right) (a-b)^2.
\label{keyineq}
\end{equation}
\end{lemma}

A direct but lengthy proof of this lemma is given in \cite{ST6}. Here, we shall give a more intuitive and rigorous proof using convexity. 

\begin{proof}  
First, let us show that
 \begin{equation}h(x):=dx^{d+\alpha}-(d+\alpha)x^d+\alpha-\alpha(2x+1)(x-1)^2\ge 0.\label{key1}
 \end{equation}
for $x\ge 0,$ either $d=2$ and $1\le\alpha\le 2$ or $3\le d\in \mathbb{N}$ and $0<\alpha\le 2.$

 \noindent\textbf{Case 1:} Assume that $d\geq 3$ and $0< \alpha \leq 2$.
 Observe that $h$ can be written as $h(x)=x^2g(x)$ where
 $$g(x)=dx^{d+\alpha-2}-(d+\alpha)x^{d-2}-2\alpha x+3\alpha.$$
 Differentiating we get
 $$g'(x)=d(d+\alpha-2)x^{d+\alpha-3}-(d+\alpha)(d-2)x^{d-3}-2\alpha,$$
 \begin{equation*}
 g''(x) = x^{d-4}d(d+\alpha-2)(d+\alpha-3) \left(x^{\alpha}-\frac{(d+\alpha)(d-2)(d-3)}{d(d+\alpha-2)(d+\alpha-3)}\right).
 \end{equation*}
 Note that $g''(x_{d,\alpha})=0$ where
 $$x_{d,\alpha}:=\left(\frac{(d+\alpha)(d-2)(d-3)}{d(d+\alpha-2)(d+\alpha-3)}\right)^{1/\alpha} < 1.$$
 When $x\ge x_{d,\alpha},$ we have $g''(x)\ge 0$, implying that 
 $g$ is convex on $[x_{d,\alpha},\infty)$. Thus, 
 $$g(x)\geq g(1)+g'(1)(x-1)=0,$$
 since $g(1)=0$ and $g'(1)=0$.
 That is, $g(x)\geq 0$ on $[x_{d,\alpha},\infty)$. In particular, $g(x_{d,\alpha})\ge 0.$ 
 On the other hand, when $0\le x\le x_{d,\alpha},$ we have
 $g''(x)\le 0$ yielding that $g'$ is decreasing on $[0, x_{d,\alpha}]$. This implies that $g'(x)\le g'(0)=-2\alpha<0,$ meaning that $g$ is decreasing on $[0, x_{d,\alpha}].$ This leads to $g(x)\ge g(x_{d,\alpha})\ge 0$ for $x\in[0, x_{d,\alpha}].$
 Hence, $g(x)\geq 0$ for $x\in[0,\infty)$. Therefore, we deduce that $h(x)=x^2g(x)\geq 0$ for $x\in[0,\infty).$ 
 \begin{figure}[!htbp]
 \centering
 \includegraphics[scale=0.6]{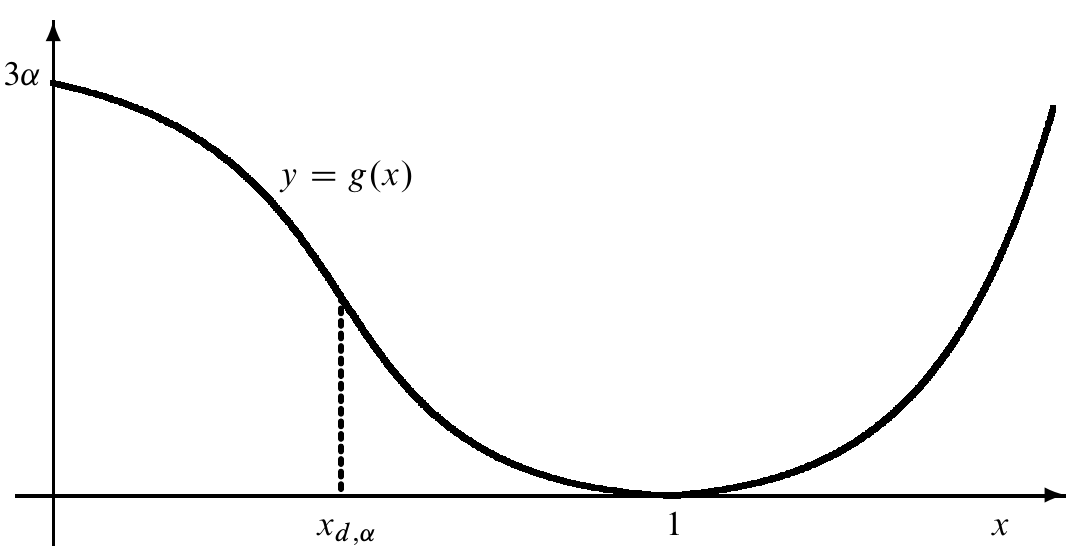}\hspace{0.2cm}
 \includegraphics[scale=0.6]{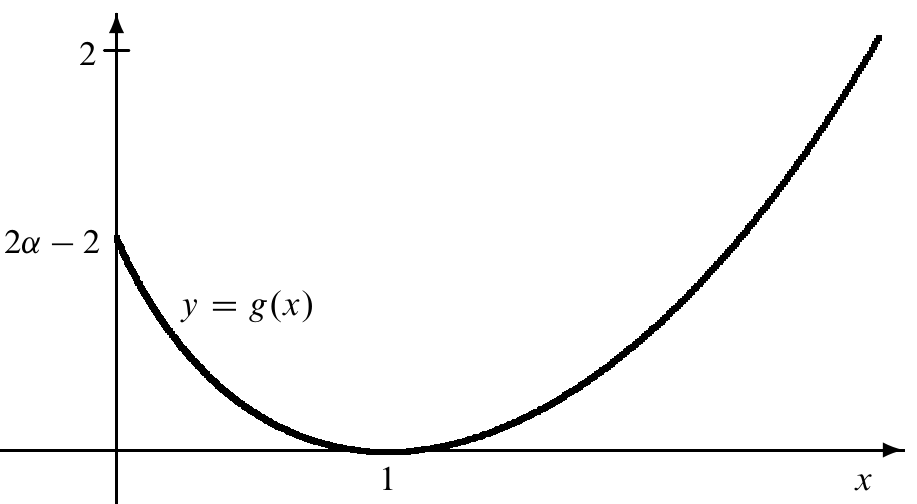}
 \caption{Graphs of $g(x)$ for $d\ge 3$ and $d=2$, 
 respectively.}
 \end{figure}
 
 \noindent\textbf{Case 2:} Now, assume that $d=2$ and $1\leq \alpha\leq 2$. Then $h$ becomes 
 $$h(x)=2x^{2+\alpha}-(2+\alpha)x^2+\alpha-\alpha(2x+1)(x-1)^2.$$
 As before, we can write $h(x)=x^2g(x)$ where $g(x)=2x^{\alpha}-2\alpha x +(2\alpha-2).$ Note that if $\alpha=1,$ then $g(x)=0.$
 Differenting again, we obtain
 $$g'(x)=2\alpha x^{\alpha-1}-2\alpha\quad\mbox{and}\quad 
 g''(x)=2\alpha(\alpha-1)x^{\alpha-2}.$$
Notice that for $1\leq \alpha\leq 2$, $g''(x)\geq 0,$ meaning that $g$ is convex. Since $g(1)=0$ and $g'(1)=0$,
 $g(x)\geq g(1)+g'(1)(x-1)=0$
 implies that $g(x)\geq 0$ and, therefore, $h(x)\geq 0.$
 
Setting $x=a/b$ in \eqref{key1}, multiplying through by $b^{d+\alpha}/d$ and rearranging the terms, we conclude \eqref{keyineq}, as desired.
 \end{proof}


\begin{remark}
When $0<\alpha<1$ and $d=2,$ the inequality above fails to hold, therefore, we do not resolve this case in this manuscript.
\end{remark}

With Lemma \ref{lem:ineq} in hand, we are now ready to prove Proposition \ref{impBLY}.

\begin{proof}
Assume that \eqref{Fintk}-\eqref{gradUk} hold and $d\geq 2.$ 
Consider the decreasing, absolutely continuous function $\varphi_k :[0,\infty)\to [0,\infty)$  defined by \eqref{zetak}. We know  that
$0\leq-\varphi_k'(t)\leq \beta$ for $t\ge 0$ where $\beta>0$ is given by \eqref{gradUk}. 
Since $\varphi_k (0)>0$ due to \eqref{Uk} let us first define
\begin{equation}
\label{trk}\Theta_k(a):=\frac{1}{\varphi_k (0)}\varphi_k \left(\frac{\varphi_k (0)}{\beta}a\right).
\end{equation}
Note that $\Theta_k$ is positive, $\Theta_k(0)=1$ and $0\le -\Theta_k'(a)\le 1$. To simplify the notation, we also set $\theta_k(a):=-\Theta_k'(a)$ for $t\geq 0$.
Hence, $0\le \theta_k(a)\le 1$ for $t\ge 0$ and  $$\displaystyle{\int_{0}^{\infty}}\theta_k(a)\,da=\Theta_k(0)=1.$$ Now,
set
\begin{equation}\label{defint}
\zeta_k=\int_0^{\infty} a^{d-1}\,\Theta_k(a)\,da\qquad\mbox{and}\qquad \eta_k=\int_0^{\infty} a^{ \alpha+d -1}\,\Theta_k(a)\,da.
\end{equation}
Using \eqref{Fintk} we get
\begin{equation}\label{eqnkA}
k=\int_{{\mathbb R}^d} \Phi_k(\mu)\,d\mu
=\int_{{\mathbb R}^d} \Phi_k^*(\mu)\,d\mu
=d\omega_d\int_{0}^{\infty}a^{d-1}\varphi_k(a)\,da.\end{equation}
Moreover, since the map $\mu\mapsto|\mu|^{\alpha }$ is radial and increasing, by \eqref{Uksumlambda}, we obtain that
\begin{equation}\label{sumxiB}
\begin{split}
\sum_{j=1}^k\lambda_j^{(\alpha)}&=\int_{{\mathbb R}^d} |\mu|^{\alpha }\,\Phi_k(\mu)d\mu
\\ &\ge \int_{{\mathbb R}^d} |\mu|^{\alpha }\,\Phi_k^*(\mu)d\mu\\ &=d\omega_d\int_{0}^{\infty} a^{ \alpha+d-1} \varphi_k(a)\,da. \end{split}\end{equation}
Substitution of \eqref{trk} into \eqref{defint} yields
\begin{equation}\label{imp01}\begin{split}
\zeta_k &=
\frac{\beta^{d}}{\varphi_k(0)^{d+1}}\int_0^{\infty} a^{d-1}\varphi_k(a)\,da=\frac{\beta^{d}k}{d\,\omega_d\,\varphi_k(0)^{d+1}},\\
\eta_k&=
\frac{\beta^{ \alpha+d }}{\varphi_k(0)^{ \alpha+d +1}}\int_0^{\infty} a^{ \alpha+d -1}\varphi_k(a)\,da\le \frac{\beta^{ \alpha+d }\,\sum_{j=1}^k\lambda_j^{(\alpha)}}{d\,\omega_d\,\varphi_k(0)^{ \alpha+d +1}}
\end{split}
\end{equation} 
Observe that application of Fubini's theorem together with $$\Theta_k(b)=\int_b^{\infty} \theta_k(a)\, da$$ leads to
\begin{eqnarray*}\frac{1}{t+d}\int_0^{\infty}a^{t+d}\,\theta_k(a)\,da&=& \int_0^{\infty}\left(\int_0^a b^{t+d-1}\, db\right)\theta_k(a)\,da\\ &=& \int_0^{\infty}b^{t+d-1}\left(\int_b^{\infty} \theta_k(a)\, da\right)db\\ &=& \int_0^{\infty}b^{t+d-1} \Theta_k(b)\, db,\end{eqnarray*} which together with $y=0$ and $y=\alpha$ respectively yields 
\begin{equation}\label{crucial1}
\int_0^{\infty}a^{d}\,\theta_k(a)\,da=\zeta_k\, d
\qquad\mbox{and}\qquad 
\int_0^{\infty}a^{ \alpha+d }\,\theta_k(a)\,da=\eta_k\,( \alpha+d ).
\end{equation}
Notice that 
\begin{equation}\big(a^{d}-1\big)\big(\theta_k(a)-\mathds{1}_{[0,1]}(a)\big)\ge 0, \qquad a\in[0,\infty).
\label{ineq:00}\end{equation}
Integrating \eqref{ineq:00} from $0$ to $\infty$ gives $$ \int_0^{\infty}a^{d}\theta_k(a)\,da \ge \frac{1}{d+1}=\gamma_d(0),$$
where $\gamma_d:[0,\infty)\to(0,\infty)$ is defined by
$$\gamma_d(x)=\int_x^{x+1}a^{d}\,da.$$
Since $\gamma_d$ is continuous and non-decreasing and $\gamma_d(x)\to \infty$ as $x\to \infty,$ the Intermediate Value Theorem provides us with the existence of $\tau\ge 0$ such that $$\gamma_d(\tau)=
\int_\tau^{\tau+1}a^{d}\,da = \int_0^{\infty}a^{d}\,\theta_k(a)\,da,
$$ which, by \eqref{crucial1}, concludes that 
\begin{equation}\label{dk}
\int_{\tau}^{\tau+1}a^{d}\,da=d\,\zeta_k.
\end{equation}
Now consider the polynomial $$T(x)=x^{ \alpha+d }-\nu_1x^{d}+\nu_2=x^d(x^{\alpha }-\nu_1)+\nu_2$$
where $$\nu_1=\frac{(\tau+1)^{ \alpha+d }-\tau^{ \alpha+d }}{(\tau+1)^{d}-\tau^{d}}>0,\qquad \nu_2=\frac{(\tau+1)^{ \alpha+d }-\tau^{ \alpha+d }}{(\tau+1)^{d}-\tau^{d}}\tau^d-\tau^{ \alpha+d }\ge 0$$
are chosen so that $T(\tau)=0$ and $T(\tau+1)=0$ and $T$ remains negative on $(\tau,\tau+1)$ and positive on $[0,\infty)\backslash [\tau,\tau+1].$
It is immediate to observe that
\begin{equation}
T(a)\left(\mathds{1}_{[\tau,\tau+1]}(a)-\theta_k(a)\right)\le 0\quad\mbox{on}\quad [0,\infty).
\label{ineq:01}
\end{equation}
Integration of \eqref{ineq:01} on $[0,\infty)$ leads to
$$\int_{\tau}^{\tau+1}a^{ \alpha+d }\,da\le \int_0^{\infty} a^{ \alpha+d }\,\theta_k(a)\,da-\nu_1\left( \int_0^{\infty} a^{d}\,\theta_k(a)\,da-\int_{\tau}^{\tau+1}a^{d}\,da\right),$$
simplifying to
\begin{equation}
\int_{\tau}^{\tau+1} a^{ \alpha+d }\,da\le \int_0^{\infty} a^{ \alpha+d }\,\theta_k(a)\,da.
\label{eqn:kv}
\end{equation}
Using \eqref{crucial1}, we infer that
\begin{equation}\label{do}
\int_{\tau}^{\tau+1}a^{ \alpha+d }\,da\leq\eta_k \,( \alpha+d ).
\end{equation}
Observe that 
\begin{equation}\label{dkJ}
d\,\zeta_k=\int_{\tau}^{\tau+1}a^{d}\,da\ge \int_{0}^{1}a^d\,da= \frac{1}{d+1}.
\end{equation}
Notice that \eqref{keyineq} gives the key  inequality in the proof of this lemma. Indeed, integrating \eqref{keyineq} in $p$ from $\tau$ to $\tau+1$ we obtain
\begin{equation}\label{m1}\begin{split}
\int_{\tau}^{\tau+1}a^{ \alpha+d }\,da&\ge \frac{ \alpha+d }{d}b^{\alpha }\int_{\tau}^{\tau+1}a^{d}\,da-\frac{\alpha }{d}b^{ \alpha+d }\\ &\quad +\frac{\alpha }{d} b^{ \alpha+d -3}\int_{\tau}^{\tau+1}(2a+b)(a-b)^2\,da.
\end{split}
\end{equation}
Note that \cite{ST4}
\begin{equation*}\label{m31}
\int_{\tau}^{\tau+1}(a-b)^2\,da\ge\min_{\tau\ge 0,\; b\ge 1/2}\int_{\tau}^{\tau+1}(a-b)^2\,da=\frac{1}{12}.\end{equation*}
\begin{equation*}\label{m03}
\int_{\tau}^{\tau+1}a\,(a-b)^2\,da
\ge \min_{\tau\ge 0,\; b\ge 1/2} \int_{\tau}^{\tau+1}a\,(a-b)^2\,da
\ge \frac{1}{2}b^2-\frac{2}{3}b+\frac{1}{4}
\end{equation*}
and so, we have
\begin{equation}\label{m3}
\int_{\tau}^{\tau+1}(2a+b)\,(a-b)^2\,da\ge b^2-\frac{5}{4}b+\frac{1}{2}
\end{equation}
for any $b\ge1/2$ and $\tau\ge 0.$
Since $(\zeta_k d)^{1/d}\ge (d+1)^{-1/d}\ge 3^{-1/2}\ge 1/2$ due to \eqref{dkJ},  setting $b=(\zeta_k d)^{1/d}$ and using \eqref{dk} and \eqref{m3}, we deduce that \eqref{m1} yields to
\begin{equation}\label{m4}
\begin{split}
\eta_k &\ge \frac{1}{ \alpha+d }(\zeta_k d)^{1+\frac{\alpha }{d}}+\frac{\alpha}{d( \alpha+d )}(\zeta_k d)^{1+\frac{\alpha -1}{d}}\\ &\quad -\frac{5\alpha}{4d( \alpha+d )}(\zeta_k d)^{1+\frac{\alpha -2}{d}}+\frac{\alpha}{2d( \alpha+d )}(\zeta_k d)^{1+\frac{\alpha -3}{d}}.
\end{split}
\end{equation}
Equations in \eqref{imp01} combined with \eqref{eqnkA} and \eqref{sumxiB} turn \eqref{m4} into
\begin{equation}\label{zeta}
\begin{split}
\sum_{j=1}^{k}\lambda_j^{(\alpha)}&\geq \frac{d}{ \alpha+d }\,\omega_d^{-\frac{\alpha }{d}} \,\varphi_k(0)^{-\frac{\alpha }{d}}
\,k^{1+\frac{\alpha }{d}}
\\
&\quad +\frac{\alpha}{( \alpha+d )}\beta^{-1}\,\omega_d^{-\frac{\alpha -1}{d}} \,\varphi_k(0)^{1-\frac{\alpha -1}{d}}
\,k^{1+\frac{\alpha -1}{d}} \\
&\quad -\frac{5\alpha}{4(\alpha+d )}\beta^{-2}\,\omega_d^{-\frac{\alpha -2}{d}} \,\varphi_k(0)^{2-\frac{\alpha -2}{d}}
\,k^{1+\frac{\alpha -2}{d}}\\
&\quad +\frac{\alpha}{2(\alpha+d )}\beta^{-3}\,\omega_d^{-\frac{\alpha -3}{d}} \,\varphi_k(0)^{3-\frac{\alpha -3}{d}}
\,k^{1+\frac{\alpha -3}{d}}.
\end{split}
\end{equation}
To finish the proof of Proposition \ref{impBLY} we shall minimize the right side of \eqref{zeta} over $\varphi_k(0).$  
 To prove \eqref{imbound1}, we show that the function, which we call $\vartheta(x)$, on the right-hand side of \eqref{zeta} with $x:=\varphi_k(0)>0$ decreases on $(0,\Omega_d].$ By \eqref{fksiless} we know that $0< x\le \Omega_d.$
 First, define 
 \begin{equation*}\label{Sen}
 \vartheta(x)= \vartheta_1(x)+\vartheta_2(x)
 \end{equation*}
 where
\begin{equation}\label{psi1}
\vartheta_1(x)=\frac{d\,
k^{1+\frac{\alpha }{d}}}{( \alpha+d )\,\omega_d^{\frac{\alpha }{d}}} x^{-\frac{\alpha}{d}}
 +\frac{\alpha 
 \,k^{1+\frac{\alpha -1}{d}}}{( \alpha+d )\beta\,\omega_d^{\frac{\alpha -1}{d}}} x^{1-\frac{\alpha-1}{d}}
 \end{equation}
and
\begin{equation}\label{psi2}
\vartheta_2(x)=\frac{\alpha
\,k^{1+\frac{\alpha -3}{d}}}{2( \alpha+d )\beta^{3}\,\omega_d^{\frac{\alpha -3}{d}}} x^{3-\frac{\alpha-3}{d}}
 -\frac{5\alpha 
 \,k^{1+\frac{\alpha -2}{d}}}{4(\alpha+d )\beta^2\,\omega_d^{\frac{\alpha -2}{d}}} x^{2-\frac{\alpha-2}{d}}
\end{equation}
Thus, it is enough to show that $\vartheta_1$ and $\vartheta_2$ defined by  \eqref{psi1} and \eqref{psi2} are also decreasing on the interval $(0,\Omega_d]$. Differentiating $\vartheta_1$ and $\vartheta_2,$ we observe that $\vartheta_1(x)$ and $\vartheta_2(x)$ are decreasing on the interval $(0,x_1)$ and $(0,x_2),$ respectively, where 
$$x_1=\left(\frac{d\,  \beta \,k^{\frac{1}{d}}}{2(d-\alpha+1)\,
\omega_d^{\frac{1}{d}}}\right)^{\frac{d}{d+1}},\qquad x_2= \left(\frac{5(2d+2-\alpha)\,  \beta \,k^{\frac{1}{d}}}{2(3d+3-\alpha)\,
\omega_d^{\frac{1}{d}}}\right)^{\frac{d}{d+1}}.$$ Hence, we particularly obtain that  
$\vartheta$ is decreasing on $(0,\Omega_d]$ when we have 
$\Omega_d \le
\min\left\{x_1,x_2 
\right\}
$
for any $k\ge 1.$ Since $x\mapsto \Gamma(x)$ is increasing for $x\ge 2$ we obtain that
\begin{equation}\Gamma\left(1+\frac{d}{2}\right)\ge \Gamma(2)=1.\label{gamsiz}\end{equation} 
In view of \eqref{gamsiz}, $\beta\ge  (2\pi)^{-d}\,\omega_d^{-\frac{1}{d}}|\Omega|^{\frac{d+1}{d}}$ and definition of $\omega_d$ given in \eqref{wd} we observe that
$$x_1\ge \left(\frac{d\, |\Omega|^{\frac{d+1}{d}} \,k^{ \frac{1}{d}}}{(d-\alpha+1)\,(2\pi)^{d}\,
\omega_d^{\frac{2}{d}}}\right)^{\frac{d}{d+1}} \ge\left(\frac{\left[\Gamma\left(1+\frac{d}{2}\right)\right]^{\frac{2}{d}} \,|\Omega|^{\frac{d+1}{d}} k^{\frac{1}{d}}}{(2\pi)^{d+1}
}\right)^{\frac{d}{d+1}}\ge \Omega_d,$$
as $k\ge 1$ and $d-\alpha+1\le 2d.$ Similarly, we obtain that
$$x_2\ge \left(\frac{5(2d+2-\alpha)\,|\Omega|^{\frac{d+1}{d}}   \,k^{\frac{1}{d}}}{2(3d+3-\alpha)\,(2\pi)^{d}\,
\omega_d^{\frac{2}{d}}}\right)^{\frac{d}{d+1}}\ge
\left(\frac{\left[\Gamma\left(1+\frac{d}{2}\right)\right]^{\frac{2}{d}}\,|\Omega|^{\frac{d+1}{d}}   \,k^{\frac{1}{d}}}{(2\pi)^{d+1}\,
}\right)^{\frac{d}{d+1}}\ge \Omega_d,$$
as $k\ge 1$ and $5(2d+2-\alpha)\ge (3d+3-\alpha).$ In conclusion, we obtain that $\vartheta(x)\ge \vartheta(\Omega_d)$ as $\vartheta$ is decreasing on $(0,\Omega_d]$. Substitution of $\beta=2(2\pi)^{-d}\,\sqrt{|\Omega|\,{\mathcal I}(\Omega)}$
given in \eqref{gradUk} together with $\varphi_k(0)=\Omega_d$
turns  \eqref{zeta} into \eqref{imbound1}. 
\end{proof}
\begin{remark}\label{rmkgen}
In view of the recent work \cite{Kim,SV}, it is worth noting that one can easily extend this for elliptic operators ${\mathcal E}_f$ defined by a kernel $f$
\begin{equation}\label{elliptic}
{\mathcal E}_f u(x)=\lim_{\epsilon\to0^+}\int_{\{|y|>\epsilon\}}
(u(x+y)-u(x))\,f(y)\,dy,
\end{equation}
where $f$ satisfies
\begin{equation} f(y)\ge \sigma \frac{{\mathcal A}_{d,\alpha}}{|y|^{d+\alpha}},\label{ker}\end{equation} ${\mathcal A}_{d,\alpha}$ is the normalizing constant in the fractional Laplacian definition \eqref{fracdef1} and $\sigma>0.$ To this end, let us consider the eigenvalue problem defined by 
\begin{equation}\label{genel}
\begin{split}
-{\mathcal E}_f\,\phi_j &=\lambda_j\,\phi_j \quad \hbox{in} \,\,\Omega,\\ 
 \phi_j  &=0 \quad\hbox{in} \,\,{\mathbb R}^d\backslash \Omega 
\end{split}
\end{equation}
It is shown in \cite{SV} that the spectrum of ${\mathcal E}_f$ is also discrete and the eigenvalues $\{\lambda_j\}_{j=1}^{\infty}$ (including multiplicities) can be sorted in an increasing order. Also, the set of Fourier transforms $\{\hat{\phi}_j\}_{j=1}^{\infty}$ of $\{\phi_j\}_{j=1}^{\infty}$ forms an orthonormal set
in $L^2({\mathbb R}^d)$ since the set of eigenfunctions $\{\phi_j\}_{j=1}^{\infty}$ is an orthonormal set in $L^2(\Omega).$ Note that we use the same notation for eigenvalues and eigenfunctions to illuminate the striking similarities though they might be different for each ${\mathcal E}_f.$ Defining $\Phi_k$ as in \eqref{Uk}, we obtain \eqref{Fintk} immediately. However, \eqref{Uksumlambda} needs to be re-written as the following inequality
\begin{equation}
\label{Nsumlambda}
\begin{split}
\sum_{j=1}^k \lambda_j&=\sum_{j=1}^k \langle \phi_j,\lambda_j \phi_j\rangle  =\sum_{j=1}^k \langle \phi_j,-{\mathcal E}_f\phi_j\rangle =\sum_{j=1}^k \int_{{\mathbb R}^d} \varrho_{\alpha}(\mu)\,|\hat{\phi}_j(\mu)|^2\,d\mu \\ & \ge \sigma\int_{\mathbb{R}^d} |\mu|^{\alpha }\, \Phi_k(\mu)\, d\mu.
\end{split}
\end{equation}
where we used the fact (Proposition 3.3. in \cite{Vald}) that $$\varrho_{\alpha}(\mu)=\int_{{\mathbb R}^d} (1-\cos(y\cdot \mu))\,f(y)\,dy\ge \sigma\,{\mathcal A}_{d,\alpha}\int_{{\mathbb R}^d} \frac{1-\cos(y\cdot \mu)}{|y|^{d+\alpha}}\,dy =\sigma |\mu|^{\alpha}.$$
Having \eqref{Nsumlambda} in hand, \eqref{sumxiB} changes as follows:
\begin{equation}\label{NsumxiB}
\begin{split}
\sum_{j=1}^k\lambda_j&\ge \sigma\int_{{\mathbb R}^d} |\mu|^{\alpha }\,\Phi_k(\mu)d\mu
\\ &\ge \sigma \int_{{\mathbb R}^d} |\mu|^{\alpha }\,\Phi_k^*(\mu)d\mu\\ &=\sigma d\omega_d\int_{0}^{\infty} a^{ \alpha+d-1} \varphi_k(a)\,da. \end{split}\end{equation}
and proceeding exactly as before using \eqref{NsumxiB} in place of \eqref{sumxiB}, we immediately arrive at the following remarkable estimate:
\begin{corollary}\label{CimpBLY} 
For $k\ge 1$, and either $1\le\alpha\le 2$ and $d=2$ or  $0<\alpha\le 2$ and $d\ge 3,$ the eigenvalues $\{\lambda_j\}_{j=1}^{\infty}$ of  \eqref{genel} defined on $\Omega \subset {\mathbb R}^d$ satisfy 
 \begin{equation}\label{Cimbound1}
 \begin{split}
 \sum_{j=1}^k\lambda_j
 &\ge(4\pi)^{\frac{\alpha}{2}}\,\frac{\sigma\,d}{ \alpha+d }
 \left(\frac{\Gamma\left(1+\frac{d}{2}\right)}{ |\Omega|}\right)^{\frac{\alpha }{d}}k^{1+\frac{\alpha }{d}}
 \\& \quad+\frac{\sigma\alpha}
  {2( \alpha+d)} \frac{|\Omega|^{\frac{1}{2}-\frac{\alpha-1}{d}}\,\Gamma\left(1+\frac{d}{2}\right)^{\frac{\alpha-1}{d}}}{(4\pi)^{\frac{1}{2}-\frac{\alpha}{2}}\,{\mathcal I}(\Omega)^{\frac{1}{2}}}\,k^{1+\frac{\alpha-1}{d}}
 \\& \quad - \frac{5\sigma\alpha}
  {16( \alpha+d )}\frac{|\Omega|^{1-\frac{\alpha-2}{d}}\,\Gamma\left(1+\frac{d}{2}\right)^{\frac{\alpha-2}{d}}}{(4\pi)^{1-\frac{\alpha}{2}}\,{\mathcal I}(\Omega)}\,k^{1+\frac{\alpha-2}{d}}
  \\& \quad+\frac{\sigma\alpha}{16( \alpha+d )}
 \frac{|\Omega|^{\frac{3}{2}-\frac{\alpha -3}{d}}\,\Gamma\left(1+\frac{d}{2}\right)^{\frac{\alpha -3}{d}}}{ {(4\pi)^{\frac{3}{2}-\frac{\alpha}{2}}\,\mathcal I}(\Omega)^{\frac{3}{2}}}\,k^{1+\frac{\alpha -3}{d}}. 
 \end{split}
 \end{equation}
\end{corollary}
Note that this estimate also improves the main result of  \cite{Kim}, which is simply the multiple of the lower bound stated in \eqref{ccmbiz} \cite{ST2} by $\sigma$.

\end{remark}


\end{document}